\documentclass[11pt]{amsart}
\usepackage{amssymb, latexsym}
\newtheorem{defin}{}
\newtheorem{saetze}[defin]{}
\newtheorem{conjec}[defin]{}
\newtheorem{lemmas}[defin]{}
\newtheorem{folger}[defin]{}
\newtheorem{bemerk}[defin]{}

\begin{document}
\title[Decision and Search in Non-abelian Cramer Shoup Public Key Cryptosystem]
{Decision and Search in Non-abelian Cramer Shoup Public Key Cryptosystem}
\author[Delaram Kahrobaei]{Delaram Kahrobaei}
\address{Mathematics Department, New York City College of Technology (CUNY), 300 Jay Street, Brooklyn, NY 11201\\
Doctoral Program in Computer Science, CUNY Graduate Center, 365 Fifth Avenue, New York, NY 10016}
\email{DKahrobaei@gc.Cuny.edu}
\urladdr{https://wfs.gc.cuny.edu/DKahrobaei/www/}
\author[Michael Anshel]{Michael Anshel}
\address{Department of Computer Science, City College of New York, 138th Street and Convent Ave, New York, New York 10031 \\
Doctoral Program in Computer Science, CUNY Graduate Center, 365 Fifth Avenue, New York, NY 10016}
\email{csmma@cs.ccny.cuny.edu }
\urladdr{http://www-cs.ccny.cuny.edu/~csmma/}

\date{July 2009}

\thanks{The first author was supported by a grant from research foundation of CUNY (PSC-CUNY) and the Faculty Fellowship Publication Program Award from City University of New York.}

\maketitle
\begin{abstract}
A method for non-abelian Cramer-Shoup cryptosystem is presented. The
role of decision and search is explored, and the platform of
solvable/polycyclic group is suggested. In the process we review
recent progress in non-abelian cryptography and post some open problems that naturally arise from this path of research.
\end{abstract}
\section{Introduction}
The field of combinatorial group theory began with decision problems of Max Dehn from 1912, known as
the \emph{word problem}, the \emph{conjugacy problem} and the \emph{isomorphism problem}. These fields have developed close connections to topology, logic and computer science.  Word problem: Let $G$ be a group given by a finite presentation. Does there exist an algorithm to determine if an arbitrary word $w$ in the generators of $G$ whether or not $w =_G 1?$ Conjugacy problem: Let $G$ be a group given by a finite presentation. Does there exist an algorithm
to determine if an arbitrary pair of words $u$ and $v$ in the generators of $G$ whether or not $u$ and $v$ define conjugate elements of $G$? By the celebrated theorem of Novikov \cite{N55} and Boone \cite{B59}, there are groups for which these questions are undecidable: they cannot be answered algorithmically. Nevertheless, because of the practical importance of these problems, a lot of effort is devoted to the development of methods for investigating finitely presented groups. The study of designing such algorithms and implementing them, is computational group theory.
These topics have generated much attention and proved of importance in modern cryptography, namely algebraic cryptography, initiated in 1999 by I.Anshel, M.Anshel and D.Goldfeld \cite{AAG99}. Algebraic key establishment protocols based on the difficulty of solving equations over algebraic structures are described as a theoretical basis for constructing public-key cryptosystems. These applications rely strongly on complexity of decision problems in
combinatorial group theory. Needless to say, complexity problems are some of the most important problems in mathematics.

\subsection{Motivation}
Key exchange problems are of central interest in cryptology
\cite{MSU}. The basic aim of key exchange problems is that two
people who can only communicate via an insecure channel want to find
a common secret key. Key exchange methods are usually based on
one-way functions; that is, functions which are easy to compute,
while their inverses are difficult to determine (see \cite{JR05,
OG08}). Here 'easy' and 'difficult' can mean that the complexities
or the practicality of the methods are far away from each other;
ideally, the one-way function has a polynomial complexity and its
inverse has an exponential complexity. Many of the known one-way
functions have a common problem: it is often easy to find a one-way
function with a polynomial complexity, but showing that there is no
inverse function with similar complexity or practicality is usually
the difficult part of the project, since the best inverse function
might just not have been discovered yet. Hence it is of interest to
investigate new one-way functions. The decision problems in
combinatorial group theory have shown much  potential for this
purpose. Another reason is that, using Shor�s algorithm the discrete
log problem and prime factorization problem admit polynomial-time
quantum algorithm. That leaves the current cryptosystem in danger if
the quantum computers are to be built! Non-abelian (a.k.a
Non-commutative) group theorists have been working in the field of
non-commutative cryptography for about ten years, but a class of
groups which provides a provably secure basis for the
non-commutative protocols is not known yet. The groups of the
authors interest are mainly polycyclic groups. These are natural
generalizations of cyclic groups but are much more complex in their
structure. Hence their algorithmic theory is more difficult.\\

In the last decade the Braid group cryptography has been of great interest to many researchers in the field and has been investigated by private and public sectors. No such investigation has been done for solvable/polycyclic groups.\\

The reason for developing a non-abelian Cramer-Shoup is to strengthen the non-commutative ElGamal \cite{khan}. We think that it should be secured against CCA1 and CCA2. The tested modification would be tested while such research can go on for developing the non-commutative ElGamal \cite{khan}.
It has been pointed oout to us that the goal of Cramer-Shoup was to defend ElGamal against chosen ciphertext attacks. We do not know if the non-commutative ElGamal \cite{khan} can be attacked using chosen ciphertext attack. However in the meantime a C-S like cryptosystem can be developed to defend the KK06.
\section{Applications of Non-commutative Group Theory in Cryptography}
In 1984, Wagner et al. \cite{wagner} proposed an approach to design
of public-key cryptosystems based on groups and semigroups with
undecidable word problem. In 2005, Birget et al. \cite{birget}
pointed out that Wagner's idea is actually not based on word
problem, but on, generally easier, premise problem. And finally,
Birget et al. proposed a new public-key cryptosystem which is based
on finitely presented groups with hard word problem. In 1999, Anshel
et al. \cite{AAG99} proposed a compact algebraic key establishment
protocol. The foundation of their method lies in the difficulty of
solving equations over algebraic structures, in particular
non-commutative groups. In their pioneering paper, they also
suggested that braid groups may be good alternative platforms for
PKC(public-key cryptography). Subsequently, Ko et al. \cite{KoLee} and Dehornoy \cite{PD04} developed the theory and practice of braid-based cryptography. The security
foundation is that the conjugator search problem (CSP) is
intractable when the system parameters, such as braid index and the
canonical length of the working braids, are selected properly. In
2002, certain homomorphic cryptosystems were constructed for the
first time for non-abelian groups due to Grigoriev and Ponomarenko
\cite{grig}. Shortly afterwards, they \cite{grigpon} extended their
method to arbitrary nonidentity finite groups based on the
difficulty of the membership problem for groups of integer matrices.
In 2004, Eick and Kahrobaei \cite{EK04} proposed a new cryptosystem
based on polycyclic groups. In 2005 Baumslag, Fine and Xu proposed
public key cryptosystem using the modular group \cite{fine},
\cite{Xu}. In 2005, Shpilrain and Ushakov \cite{SU05} suggested that
R. Thompson's group may be a good platform for constructing
public-key cryptosystems. In their contribution, the key assumption
is the intractability of the decomposition problem, which is more
general than the conjugator search problem, defined over R.
Thompson's group, also a infinite non-abelian group given by finite
presentation.  In \cite{CDW07}, Cao et al. propose a new method for
designing public key cryptosystems based on general non-commutative
rings. The key idea of their proposal is that for a given
non-commutative ring, they can define polynomials and take them as
the underlying work structure. In 2006 Kahrobaei and Khan, proposed
a non-commutative key-exchange scheme which generalizes the
classical ElGamal Cipher \cite{khan}. This scheme is closer to the spirit of ElGamal and the they proposed polycyclic groups for such protocol.

Recently, Gilman, et al in \cite{GMMU} study the algorithmic
security of the Anshel-Anshel-Goldfeld (AAG) key exchange scheme and
show that contrary to prevalent opinion, the computational hardness
of AAG depends on the structure of the chosen subgroups, rather than
on the conjugacy problem of the ambient braid group. Proper choice
of these subgroups produces a key exchange scheme which is resistant
to all known attacks on AAG.

In recent work \cite{GKN}, Grigoriev et al,  show new constructions
of cryptosystems based on group invariants and suggest methods to
make such cryptosystems secure in practice. In their paper their
introduce a new notion of cryptographic security, a provable break,
and prove that cryptosystems based on matrix group invariants and
also a variation of the Anshel-Anshel-Goldfeld key agreement
protocol for modular groups are secure against provable worst-case
break unless $NP \subset RP$.

\subsection{Non-Commutative Key Exchange using Conjugacy}
The following is the key exchange problem proposed by Kahrobaei and Khan \cite{khan}.\\
Let $G$ be a finitely presented non-abelian group having solvable
word problem.  Let $S,T < G$ be finitely generated proper
subgroups of $G$, for which the subgroup $[S,T]$ (i.e. the
subgroup generated by $ \{[s,t] \;|\; s\in S, t \in T\}$) is the
trivial subgroup consisting of just the identity element of $G$.
Now suppose two parties, Alice and Bob, wish to establish a
session key over an unsecured network.

Bob takes $s \in S, b \in G$ and publishes $b$ and $c=b^s$ as his
public keys, keeping $s$ as his private key. If Alice wishes to
send $x \in G$ as a session key to Bob, she first chooses a random
$t \in T$ and sends $$E=x^{(c^t)}$$ to Bob, along with the header
$$h=b^t.$$ Bob then calculates $(b^t)^s=(b^s)^t=c^t$ with the header.
He can now compute $$E'=(c^t)^{-1}$$ which allows him to decrypt
the session key,
$$(x^{(c^t)})^{E'}=(x^{(c^t)}){^{(c^t)}}^{-1}=x.$$The element $x
\in G$ can now be used as a session key.

The feasibility of this scheme rests on the assumption that
products and inverses of elements of $G$ can be computed
efficiently.  To deduce Bob's private key from public information
would require solving the equation $c=b^s$ for $s$, given the
public values $b$ and $c$. This is called the {\em conjugacy
search problem} for $G$. Thus the security of this scheme rests on
the assumption that there is no fast algorithm for solving the
conjugacy search problem for the group $G$.

\subsection{Non-Commutative Key Exchange using Power Conjugacy}
The following is the key-exchange problem which was proposed by Kahrobaei and Khan in \cite{khan}
What if the conjugacy search problem is tractable?  The next
paradigm embellishes conjugacy-based key exchange to address this
possibility. Bob takes $s \in S$, $g \in G$ and $n \in \mathbb{N}$
and publishes $v=g^{n}$ and $w=s^{-1}gs$ as his public keys. Note that the centralizer of $g$ is trivial. Bob
keeps $n \in \mathbb{N}$ and $s \in S$ as his private keys.  Note
that $v$ and $w$ satisfy $w^{n}=s^{-1}vs$. If Alice wishes to send
$x \in G$ to Bob, she first chooses a random $m \in \mathbb{N}$
and $t \in T$. To encrypt $x$, Alice computes
$$E=x^{-1}t^{-1}(v)^{m}tx=x^{-1}t^{-1}g^{mn}tx$$ and sends it to Bob
along with the header $$h=t^{-1}w^{m}t=t^{-1}s^{-1}g^{m}st.$$ Bob
receives $E$ and $h$, and computes $$E' =
sh^{n}s^{-1}=t^{-1}g^{mn}t.$$  Note that $E = x^{-1}E' x$, so if
Bob can solve the conjugacy search problem, he can obtain $x \in
G$, which can then serve as the common secret that can be used as
a symmetric session key for secure communication.

The feasibility of this scheme rests on the assumption that
products and inverses of elements of $G$ can be computed
efficiently, and that the conjugacy problem is solvable. To deduce
Bob's private key from public information would require solving
the equation $w^{n}=s^{-1}g^{n}s$ for $n$ and $s$, given the
public values $g^{n}$ and $w$. This is called the {\em power
conjugacy search problem} for $G$. Thus the security of this
scheme rests on the assumption that there is no fast algorithm for
solving the power conjugacy search problem for the group $G$.

\section{Cramer-Shoup cryptosystem}
Cramer-Shoup cryptosystem is a generalization of ElGamal Key exchange problems, it is provably secure against adaptive chosen ciphertext attack. Moreover, the proof of security relies only on a standard intractability assumption, namely, the hardness of the Diffie-Hellman decision problem in underlying group (see \cite{shoup}, \cite{shoup2}). A hash function $H$ whose output can be interpreted as a number in $\mathbb Z_q$ (where $q$ is a large prime number). It should be hard to find collisions in $H$. In fact, with a fairly minor increase in cost and complexity, one can eliminate $H$ altogether.
\subsection{Definition of provably secure against adaptive chosen ciphertext attack}
The right formal, mathematical definition of security against active
attacks evolved in a sequence of papers by Naor and Yung, Rackoff
and Simon, Dolev, Dwork and Naor. The notion is called chosen ciphertext security or
equivalently non-malleability.
The intuitive thrust of this definition is that even if an
adversary can get arbitrary ciphertexts of his choice decrypted, he
still gets no partial information about other encrypted messages. (see for more information see \cite{shoup}, \cite{shoup2})
\subsection{The Cramer-Shoup Scheme}
\textbf{Secret Key:} random $x_1, x_2, y_1, y_2, z \in \mathbb Z_q$\\
\textbf{Public Key:}
\begin{center}
$g_1, g_2$ in $G$ (but not $1$)\\
$c = {g_1}^{x_1}{g_2}^{x_2}, d= {g_1}^{y_1}{g_2}^{y_2}$\\
$h={g_1}^z$.\\
\end{center}
\textbf{Encryption} of $m \in G$: $(u_1, u_2, e, v)$, where
\begin{center}
$u_1 = {g_1}^r, u_2 = {g_2}^r, e = h^r m, v=c^r d^{r \alpha} , r$ in $\mathbb Z_q$ is random, and\\
$\alpha = H(u_1, u_2, e).$
\end{center}
\textbf{Decryption of $(u_1, u_2, e, v)$:}
\begin{center}
If $v= {u_1}^{x_1+ \alpha y_1} {u_2}^{x_2 + \alpha y_2}$, where $\alpha = H(u_1, u_2, e)$\\
then $m= e/{{u_1}^z}$\\
else "reject"
\end{center}

\section{Non-commutative Cramer-Shoup Key-exchange problem}
The objective of this section is to propose non-commutative Cramer-Shoup cryptosystem and analyze its security. As it is pointed out in \cite{shoup}, \cite{shoup2}, with a fairly minor increase in cost and complexity, one can eliminate the hash function altogether and that is what we are doing in this scheme. However the following problem still remains open: Can our protocol be extended that our C-S like cryptosystem can be modified which needs a Hash function?

Let $G$ be a non-abelian group such that every element has a normal form, and the conjugacy search problem is hard.\\
\textbf{Secret Key:} random $x_1, x_2, y_1, y_2, z \in G$\\
\textbf{Public Key:}
\begin{center}
$g_1, g_2$ in $G$ (but not $1$), such that $[{g_2}^{x_2},{g_1}^{y_1}]=1$\\
$c = {g_1}^{x_1}{g_2}^{x_2}$, where ${g_1}^{x_1}= {x_1}^{-1} {g_1} {x_1}$ and ${g_2}^{x_2}= {x_2}^{-1} {g_2} {x_2}$\\
$d= {g_1}^{y_1}{g_2}^{y_2}$ , where ${g_1}^{y_1}= {y_1}^{-1} {g_1} {y_1}$ and ${g_2}^{y_2}= {y_2}^{-1} {g_2} {y_2}$\\
$h={g_1}^z =z^{-1} {g_1} z$.\\
\end{center}
\textbf{Encryption} of $m \in G$: $(u_1, u_2, e, v)$, where
\begin{center}
$u_1 = {g_1}^r = r^{-1} {g_1} r, u_2 = {g_2}^r = r^{-1}{g_2}r, e = m^{h^r}, v=c^r d^r$.
\end{center}
Suppose $r$ be a random element in $G$ such it commutes with $x_1, x_2, y_1, y_2$ and $z$.\\
\textbf{Decryption of $(u_1, u_2, e, v)$:}
\begin{eqnarray*}
                                           \text{If     }      v &=& {u_1}^{x_1}{u_1}^{y_1} {u_2}^{x_2}{u_2}^{y_2} \\
{({{g_1}}^{x_1} {{g_2}}^{x_2})}^r {({g_1}^{y_1} {g_2}^{y_2})}^r &=& {{g_1}^r}^{x_1} {{g_1}^r}^{y_1} {{g_2}^r}^{x_2}{{g_2}^r}^{y_2}\\
{{g_1}^{x_1}}^r {{g_2}^{x_2}}^r {{g_1}^{y_1}}^r {{g_2}^{y_2}}^r &=& {{g_1}^{x_1}}^r {{g_2}^{x_2}}^r {{g_1}^{y_1}}^r {{g_2}^{y_2}}^r
\end{eqnarray*}
\begin{center}
then $m= e^{{({u_1}^z)}^{-1}}= m^{{{h^r}}{{{u_1}^{{z}^{-1}}}}} = m^{{{{u_1}^z}}{{({u_1}^{z})}^{-1}}}$\\
else "reject"
\end{center}
Note that the authors currently experimenting the search for the secure randomly generated $r$.


\section{Polycyclic groups new platform for cryptology}
Polycyclic groups are a natural generalization of cyclic groups, but
they are much more complex in their structure than cyclic groups.
Hence their algorithmic theory is more difficult and thus it seems
promising to investigate classes of polycyclic groups as candidates
to have a more substantial platform perhaps more secure. Recall that
a group is called polycyclic if there exists a polycyclic series
through the group; that is, a subnormal series of finite length with
cyclic factors. There are two different natural representations for
these groups which can be used for computations: polycyclic
presentations and matrix groups over the integers. We refer to
\cite{S94}, \cite{segal}, \cite{poly}, \cite{BE01}, \cite{BCRS91} and \cite{holt} for
background and a more detailed introduction to polycyclic groups.

In particular polycyclic groups are linear. In this setting, both group multiplication and the word
problem are efficiently solvable, since matrix multiplication for such groups is
solvable in polynomial time.
One can show that for a subgroup of a general linear group, if two elements are conjugate then they have the same Jordan normal form. Using this lemma one conclude that the search conjugacy problem in any subgroup
of the General Linear group is solvable. However the complexity is not known but is conjectured to be exponential.
Kharlampovich \cite{OK81} showed that there is a finitely presented solvable group with an undecidable word problem. It follows by a theorem of Arzhantseva-Osin in \cite{AO2002} that the word problem is in NP for any finitely generated metabelian group. For polycyclic groups, some decision
problems are known to be difficult but not provably so \cite{eick, LGS98}.\\
A large growth rate would imply a large key space for the set
of all possible keys, thus making the brute force search of this
space intractable.  Ideally, we would like to use groups which
exhibit provably exponential growth. A large class of
polycyclic groups are known to have an exponential growth rate,
namely those which are not virtually nilpotent, by results of Wolf
and Milnor in 1968.
Using these observations about the polycyclic groups, we conjecture they are the best candidate for the non-abelian Cramer-Shoup cryptosystem.
\bibliographystyle{amsplain}
\bibliography{XBib}
\end{document}